

\input amstex 
\documentstyle{amsppt} 
\magnification = \magstep1
\hsize = 6.25 truein
\vsize = 22 truecm
\baselineskip .22in
\define\Lap{\varDelta}
\define\del{\partial}
\define\a{\alpha}
\predefine\b{\barunder}
\redefine\b{\beta}

\predefine\d{\dotunder}
\redefine\d{\delta}

\define\e{\epsilon} 
\define\s{\sigma}
\predefine\o{\orsted}
\redefine\o{\omega}
\define\Sig{\Sigma}
\define\Lam{\Lambda}
\define\lam{\lambda}

\define\RR{\Bbb R}
\define\SS{\Bbb S}
\define\HH{\Bbb H}

\define\ML{\Cal M_{\Lambda}}

\define\MMk{\Bbb M_{k}}
\define\CML{{\overline{\Cal M}}_{\Lambda}}

\define\Min{\text{Min}}
\define\Cmc{\text{CMC}}
\define\Csc{\text{CPSC\,}}

\NoRunningHeads 
\topmatter 
\title Gluing and moduli for noncompact geometric problems \endtitle 
\author Rafe Mazzeo${}^{(\dag)}$ and Daniel Pollack${}^{(\star)}$ \endauthor  
\affil Stanford University ($\dag$) and 
University of Chicago ($\star$) \endaffil 
\thanks Research supported in part by ($\dag$) NSF Young Investigator 
Award and NSF grant \# DMS9303236, and ($\star$) NSF Postdoctoral 
Research Fellowship.
\endthanks
\endtopmatter
\document

\specialhead I. Introduction \endspecialhead

In this paper we survey a number of recent results concerning the 
existence and moduli spaces of solutions of various geometric problems on 
noncompact manifolds. The three problems which we discuss in detail are: 
$$
\align
\text{\bf I.} & \text{  Complete properly immersed minimal surfaces in 
$\RR^3$ with finite total} \\
& \text{  \qquad curvature.} \\
\text{\bf II.} & \text{  Complete embedded surfaces of constant mean 
curvature in $\RR^3$ with finite} \\
& \text{    \qquad topology.} \\ 
\text{\bf III.} & \text{  Complete conformal metrics of constant positive 
scalar curvature on $M^n \setminus \Lambda$,}\\
&\text{\quad where $M^n$ is a compact Riemannian manifold, $n\geq3$ and $\Lam 
\subset M$ is closed.}
\endalign
$$

The existence results we discuss for each of these
problems are ones whereby known solutions (sometimes satisfying certain
nondegeneracy hypotheses) are glued together to produce
new solutions. Although this sort of procedure is quite well-known,
there have been some recent advances on which we wish to report here.
We also discuss what has been established about the moduli spaces
of all solutions to these problems, and report on some work in progress
concerning global aspects of these moduli spaces.  In the final 
section we present a new compactness result for the `unmarked moduli
spaces' for problem {\bf III}.

Although the analysis underlying each of these problems differs somewhat 
from one case to the next, there are many common themes.
The basic point which makes each of these problems tractable
is the fact that the geometry of the ends of the manifolds
of interest is well-understood. 
The ends of a surface $\Sigma$ in case {\bf I}
are either asymptotically planes or catenoids, in case {\bf II} 
the ends are asymptotically periodic, and in case {\bf III},
when $\Lam$ is regular, the ends 
are either asymptotically periodic or asymptotically of `edge type',
which means that near infinity they are locally asymptotic to a
neighbourhood of infinity in a product of hyperbolic space and
a compact manifold. In each of these cases there is a well-developed
set of techniques to study the geometric linear elliptic operators
on each of these types of manifolds. 

The point of view we wish to promulgate here is that gluing
methods for constructing new solutions from old ones, and the study
of the local structure of the moduli space of all solutions of any one 
of these problems are most effectively handled by using linear
analysis tailored specifically to the geometry at hand.  In what follows
we  survey recent work in these areas, in some of which we have played
a role, and attempt to give some idea of the techniques 
we have found to be successful.  One reason for including problem {\bf I} 
here is the great similarity between the moduli space results obtained
by P\' erez and Ros \cite{PR}, using rather different methods, 
and our own results in these directions for problems {\bf II} and {\bf III}.  

In an attempt to discuss a large number of results by many authors in 
a brief survey we have undoubtedly given some topics less attention 
than they would deserve in a more thorough survey of the field.
Suffice it to say that our aim is to present those results which we 
have found most relevant in our attempts to understand and establish new 
results concerning gluing and moduli for these geometric problems.

\specialhead 2. Gluing \endspecialhead

A very rough outline of the gluing method is as follows. If one wishes to
construct a new solution of one of the problems {\bf I} -- {\bf III}, the 
first step is to construct an approximate solution. This is built out of 
pieces of actual solutions already at hand, which are somehow joined 
together.  This approximate solution is usually an actual solution away 
from the transition regions where the pieces are adjoined.
Suppose we are in cases {\bf I} or {\bf II}.
Then surfaces near to this approximate solution
surface, $\Sigma$, may be written as a normal graph off of $\Sigma$, i.e. as
$$
\Sigma_\phi = \{x + \phi(x) \nu(x): x \in \Sigma \} \tag 2.1
$$
where $\nu(x)$ is the unit normal to $\Sigma$ at $x$ and $\phi$
is small. In case {\bf III} this normal perturbation is replaced by a 
conformal factor close to $1$: metrics near to some $g_0$ are written as
$$
g = (1 + v)^{\frac{4}{n-2}} g_0,  \tag 2.2
$$
for some function $v$ which is taken to be small in an appropriate norm.
For simplicity we call these perturbations 
`normal perturbations'. We also need to consider nearby surfaces or metrics
in a slightly broader sense. Thus we also consider normal perturbations off 
of not just $\Sigma$ or $g_0$ but also surfaces or metrics where the
geometry of the ends is slightly changed and allowed to vary slightly
within the associated asymptotic family of solutions. For example, for 
the extrinsic problems {\bf I} or {\bf II} this includes small 
independent rotations or translations of each end.
We make this more precise below in each case, but the salient fact is
that these asymptotic changes form a finite dimensional
family. Thus the set of all solutions to the geometric problem
near to $\Sigma$ or $g_0$ are contained within the set of normal perturbations
off of elements of this finite dimensional family of perturbations of
the ends.

Amongst all such surfaces or metrics one seeks one which is
a solution of the geometric problem. This may be expressed as a second order
elliptic partial differential equation in $\phi$ or $v$.  
For {\bf I} and {\bf II} 
this equation is quasilinear, while for {\bf III} it is only semilinear, 
but with a critical nonlinearity.  The gluing construction is completed 
if one can find a solution of this nonlinear equation.

There are various difficulties in implementing this method. First, 
the construction of the approximate solution can be extremely 
delicate. Second, it can be quite difficult to obtain a
solution of the nonlinear equation: frequently (and for
geometric reasons) its linearization is not surjective, and
in addition one must deal with the fact that the surface or manifold
is noncompact.  As we discuss below, both of these
problems are handled in part by introducing a supplementary
finite dimensional space $W$, the deficiency subspace,
which contains the linearizations of the changes in the
geometry of the ends mentioned above. 

In the remainder of this section we discuss the three
cases in turn, surveying the existence results that have
been proved using gluing methods and discussing 
the last case in some detail. 
 
\head Embedded Minimal surfaces with finite total curvature \endhead

Amongst the three geometric problems we have mentioned, the one
that has received by far the most attention is the first, concerning
the existence and nature of complete embedded minimal surfaces
in $\RR^3$. A surface $\Sigma \subset \RR^3$ is said to be
minimal if its mean curvature $H$, which at a point $p \in \Sigma$ is 
the average of the two principal curvatures $\kappa_1$ and $\kappa_2$ 
of $\Sigma$ at $p$, vanishes. This is equivalent to saying that the 
surface is a critical point for the area functional under compactly
supported variations.  The two simplest (complete) examples
are the plane and the catenoid. The latter is the surface of 
revolution obtained by revolving the catenary, $y = \cosh x$
around the $x$-axis. There is a classical method for generating
minimal surfaces using complex function theory, leading to
what is known as the Weierstrass representation for minimal
surfaces in $\RR^3$. Given certain holomorphic data, one immediately
writes down a minimal immersion. Unfortunately this immersion
is rarely a proper embedding, and even if it is, it may well
correspond to a surface of infinite topology.  For example, there
are many explicit singly or doubly periodic minimal surfaces.
One of the central questions concerns the existence and nature of 
complete embedded minimal surfaces for which the integral
of the Gauss curvature is finite. Of the classically known
examples, only the plane and catenoid have this property.
In the early 1980's the first new example was discovered by 
C. Costa \cite{Co} who exhibited a genus one minimal immersion
with one planar and two catenoidal embedded ends.  
Costa's surfaces was later shown to be embedded by D.~Hoffman and 
W.~Meeks~III \cite{HM}.
This breakthrough led to a huge resurgence of interest
in this area. Since that time many new families of examples
of complete embedded minimal surfaces have been discovered.
We refer the reader to two very well-written sources,
the well-known monograph by R. Osserman \cite{O} and the recent survey 
paper of D.~Hoffman and H.~Karcher \cite{HKa}, for more information.

An important tool in the study of a complete embedded minimal surface
$\Sigma$ is the Gauss map $\nu: \Sigma \rightarrow \SS^2$. The 
minimality of $\Sigma$ implies that $\nu$ is conformal; this leads
to the fundamental fact, first established by Osserman,  
that $\Sigma$ is conformally equivalent to
the complement of a finite set of points $\{p_1, \dots, p_k \}$
in a {\it compact} Riemann surface $M$. In particular, $\Sigma$
must have finite topology. It is also known that
each end of $\Sigma$ must be asymptotically equivalent to 
a plane or a catenoid. In order for these ends not to intersect,
their axes must be parallel; this is equivalent to the assertion
that the limit of $\nu$ along each end is one of two
fixed antipodal points in $\SS^2$. We normalize by demanding that these
ends are asymptotically horizontal, i.e. that the limiting value
of $\nu$ at each end is $\pm e_3$. 

As noted, the construction of new minimal surfaces has relied 
primarily on the techniques associated with the Weierstrass representation.
A more satisfactory and flexible method would be provided by gluing 
techniques, and quite recently N. Kapouleas has announced 
a dramatic new result of this type which we discuss briefly.
His result provides a method to `desingularize' a transverse 
intersection of complete embedded minimal surfaces of finite
total curvature. At present he has announced a result when the
initial configuration of intersecting surfaces is rotationally invariant,
and expects the general case to be completed soon.

The desingularization rests on the existence of a family of
simply periodic complete embedded minimal surfaces which desingularize
the intersection of two planes meeting at any angle. These surfaces
were discovered by Scherk in the 19th century via their Weierstrass
representations.  They may be visualized by replacing the line
of intersection by a periodic string of `handles'. 

Now, let $\{\Sigma_j\}$ be a finite collection of complete embedded
rotationally invariant minimal surfaces of finite total curvature
sharing a common axis of rotation. Assume that these surfaces meet 
transversely and that there are no triple intersections. This set is called
the initial configuration.  The curves of intersection of the
different $\Sigma_j$ are circles.  The intersection $\Sigma_i \cap \Sigma_j$ 
is desingularized as follows. First an approximate solution is constructed
by bending a Scherk surface around this curve and joining
its planar `wings' to $\Sigma_i$ and $\Sigma_j$ appropriately.
This step requires very delicate estimates to show that this
can be done in such a way that the error terms incurred in the
transition regions are suitably small. These error terms are then
dealt with by a careful analysis of the linearized mean curvature,
or Jacobi, operator.  Kapouleas' methods for this step are rather different 
than the ones we outline below. This then allows him to solve 
the nonlinear problem and find a new complete embedded minimal surface 
with finite total curvature as a perturbation over the initial configuration.  

Recently M.~Traizet \cite{Tr} has given a closely related minimal 
desingularization of the intersection 
of a finite collection of planes with normals lying in a common two-plane, 
satisfying certain genericity conditions.

\head Constant mean curvature surfaces \endhead

The study of complete surfaces in $\RR^3$ with constant mean curvature (CMC)
is almost as old, although it is perhaps not quite as intensively studied.
An old theme in the subject was the question of whether it is possible 
for there to exist a compact embedded or immersed surface without
boundary of constant mean curvature other than the sphere. Alexandrov
\cite{A} showed that the round sphere is the only such surface which
is embedded.
For several decades this question remained open, until the early 1980's
when Wente \cite{W1} (and subsequently \cite{W2}) constructed constant 
mean curvature immersions of a torus.  This construction has since been 
reinterpreted and amplified using completely integrable systems; this 
has led to a classification of all immersed CMC tori, cf. the work of 
Pinkall and Sterling \cite{PS},  Abresch \cite{Ab1} \cite{Ab2}, 
and Bobenko \cite{B}.

It turns out that there exist constant mean curvature immersions 
of compact surfaces without boundary of all genera greater than zero.
The existence of these higher genus examples was established
in the  landmark papers of Kapouleas \cite{K2}, \cite{K3}.
These papers established the gluing method as a valuable tool
in extrinsic geometry. In the first, certain `balanced' arrangements
of spheres and pieces of Delaunay surfaces (described below)
are joined together to obtain an immersed approximate
solution, and then perturbed to an exact solution.
In \cite{K1} he uses a similar method to construct many complete 
noncompact CMC surfaces, both immersed and embedded.
In the compact case only surfaces of genus greater than two 
may be obtained this way.
To remedy this Kapouleas found, in \cite{K3}, a more elaborate and general
construction whereby Wente tori are fused to each other
and to these other pieces. Before 
describing these results more we need to digress to discuss 
the noncompact Delaunay surfaces, pieces of which, along with
spheres and Wente tori, serve as the other building blocks
in Kapouleas' construction. 

The simplest noncompact CMC surface is the cylinder. One can ask
which other surfaces of revolution have constant mean curvature,
and phrased this way it is not difficult to discover a
one-parameter family of curves such that the corresponding
surfaces of revolution have constant mean curvature. Furthermore
these are the only complete CMC surfaces of revolution besides the sphere.
These were discovered by  C. Delaunay in 1841 \cite{D},
and are commonly known as the Delaunay surfaces. The curves
producing them are periodic and the embedded surfaces in this family
may be parametrized by the minimum distance to the axis of revolution. 
If we normalize the mean curvature to be $1$, then this minimum distance 
$\e$ varies between $0$ and $1$ and the corresponding surfaces are
denoted $D_\e$. Thus $D_{1}$ is the unit cylinder and as $\e\rightarrow 
0$ the surfaces $D_{\e}$ converge to the singular CMC surface formed by 
an infinite string of mutually tangent spheres of radius 2 arranged 
along an axis. By choosing the parametrization somewhat differently
we obtain, for $\e<0$, generating curves which cross the axis of revolution
and give rise to Delaunay surfaces which are immersed.  There is a very pretty
geometric description of this whole family of generating curves: they are the 
curves traced by one focus of a conic rolling along the axis. For example,
the cylinder is obtained as the locus of the center of a rolling
circle of radius 1; the embedded generating curves are the loci
of rolling ellipses, where, as $\e\rightarrow 0$, the eccentricity of
the ellipse tends to infinity, the singular CMC surface obtained when $\e=0$
corresponds to `rolling' the infinite eccentricity conic, i.e. 
the line segment of length 2. Finally, the nonembedded generating curves 
correspond to rolling hyperbol\ae. This is explained in detail in \cite{E}.

There are many difficulties in obtaining sufficiently good approximate 
solutions in Kapouleas' construction; this step is more
subtle than in the intrinsic case {\bf III} discussed below. 
One starts with a `labeled graph', i.e. a collection of vertices and 
edges as points and line segments in $\RR^3$, such that each edge is 
assigned a real number, 
which satisfies certain `balancing conditions' involving the geometry
of the graph as well as the labeling.
One then builds an approximate solution by connecting spheres or
pieces of Delaunay surfaces with very small necksize at each vertex or 
string of vertices, as prescribed by this labeled graph. 
The edges between vertices correspond to
catenoid-shaped necks. Great care is needed in joining these
various pieces together to obtain a good approximate solution. 
The second step is to try to perturb this configuration to an exact 
solution, and this involves a thorough understanding of the 
Jacobi operator $\Cal L$. $\Cal L$ is a second order, self-adjoint 
elliptic operator which, in cases {\bf  I} and {\bf II}, takes the simple form
$$
\Cal L = \Lap_{\Sigma} + |A_\Sigma|^2,
$$
where $|A_\Sigma|^2$ is the squared norm of the second fundamental
form of $\Sigma\hookrightarrow\RR^3$.  The main complication 
in this step is that $\Cal L$ has a number of small or
vanishing eigenvalues, caused by the degeneracy, or existence of null modes, 
of the Jacobi operator on each of the building blocks, i.e. the spheres,
Wente tori and Delaunay surfaces. This degeneracy is precisely what separates
a hard gluing theorem from some of the softer ones we describe below.

The non-invertibility of $\Cal L$ on the sphere, for example, is a geometric 
phenomenon which is not difficult to understand. Let $\Sigma_t$ be a 
one-parameter family of CMC surfaces deforming $\Sigma_0=\SS^2$ and
write $\Sigma_t$ as a normal graph off of $\Sigma_0$ of some function $\Phi_t$.
This function solves the CMC equation $N(\Phi_t) = 0$; differentiating this 
and evaluating at $t=0$ produces an element of the nullspace of $\Cal L
= \Lap_{\SS^2} + 2$. The obvious deformations comes from translations
in $\RR^3$, and from these we get a 3-dimensional nullspace of $\Cal L$ on the
sphere (or any closed CMC surface in $\RR^3$). These Jacobi fields 
are the restrictions to the sphere of the linear functions in $\RR^3$. 

To get around the problem of the degeneracy of the Jacobi operator,
Kapouleas identifies an explicit approximate nullspace, and 
by carefully adjusting his approximate solution he generates a 
`substitute kernel' which in effect allows him to work orthogonally to 
the approximate nullspace. In the noncompact case this requires 
infinitely many adjustments.  Kapouleas' work is also briefly 
discussed in \cite{S4}.

It turns out, by a theorem of Korevaar, Kusner, and Solomon \cite{KKS} 
that any complete embedded CMC surface is in many ways similar to
the surfaces that Kapouleas constructs. In particular, each
end is strongly asymptotic to an embedded Delaunay surface. The proof
of this uses an Alexandrov reflection argument, and depends rather
strongly on the embeddedness (or the slightly weaker condition
of `Alexandrov embeddedness'). In further work of Kusner and Korevaar
\cite{KK}, it is shown that there exists a graph in $\RR^3$ such that 
$\Sigma$ is contained in a neighbourhood of this graph consisting of 
the union of tubular neighborhoods of radius 6 about the edges with
balls of radius 21 about the vertices. 

Very explicit symmetric examples have been constructed by Grosse-Brauckmann
\cite{G}. For example, he constructs one-parameter families of
surfaces with $k$ ends which have a $k$-fold dihedral symmetry, and consist 
of $k$ Delaunay ends joined at a center.  On one extreme of this family they
meet at a central sphere, while at the other extreme they join at a concave 
surface which is asymptotically, as the necksizes tend to zero,
$k-$noidal c.f. \cite{G}.  
Quite recently Grosse-Brauckmann and  Kusner \cite{GK} have established 
new necessary conditions for the existence of certain CMC surfaces with
three or four ends possessing some (but not full $k$-fold dihedral) symmetry.
They also have found surfaces of this type experimentally, using a computer;
rigorous proofs of their existence would be extremely interesting.

Currently the authors and Kapouleas \cite{KaMP} are developing a gluing
procedure to establish an interior bridge principle for nondegenerate
CMC surfaces. The method used here combines Kapouleas' ideas
for constructing a one-parameter family of approximate solutions with
arbitrarily small neck as the bridge and the authors techniques,
developed in \cite{MPa1}, \cite{MPU2}, for perturbing
an approximate solution on the connected sum of two nondegenerate 
solutions to an exact solution. One main ingredient is 
a uniform analysis of $\Cal L$ as the size of the bridging neck shrinks.
Unlike the scalar curvature problem discussed below in 
detail, the introduction of a small neck automatically leads to 
the existence of small eigenvalues for $\Cal L$. This is because the neck
is shaped like a catenoid asymptotically, hence the solution of 
$\Cal L \phi = 0$ on this region may be arbitrarily well approximated by 
solutions of this equation on the catenoid. But there are solutions
of this equation on the catenoid which decay at infinity; they are
produced from translations perpendicular to the axis of the catenoid
in the manner described above. 
Thus it seems impossible to develop a gluing method for CMC surfaces 
which deals only with nondegenerate surfaces, even if the original surfaces
to be glued are nondegenerate. 

The surfaces produced by this last gluing method are quite different
than the ones obtained by Kapouleas earlier. In particular, our method
produces a $(3(2k)-6)$-parameter family (up to ambient Euclidean motions) 
of CMC embeddings of the connected sum of $k$ different
Delaunay surfaces.  This produces an open set in the moduli space
of embedded CMC surfaces with $2k-$ends, even though the resulting
surfaces may be degenerate for the reasons indicated above.
The necksizes of these surfaces need not be
small, although the necks joining the different Delaunay surfaces
are quite small. The final perturbation, from the approximate
solution to an exact solution, produces a small change in the
necksizes on half of each summand, but the necksizes on the
other ends are not changed. 

Thus for every $k \ge 2$ and a 
collection of necksizes $(\e_1, \dots, \e_k)$ we prove the 
existence of a complete embedded surface with constant mean curvature 
and $2k$ ends such that $k$ of the ends have prescribed necksizes
$\e_1, \dots, \e_k$. Included amongst these are the first
examples of complete embedded CMC surfaces with cylindrical
ends which are not cylinders. Of course, only half the ends of
these surfaces are cylindrical; on the other ends a slight
periodicity has been introduced by the final perturbation.
A related open question was recently posed to us by R.~Schoen:
is a complete embedded CMC surface, all of whose ends are 
asymptotically  cylindrical (or equivalently, having finite total
absolute Gauss curvature) necessarily a cylinder?

\head Complete metrics of constant positive scalar curvature \endhead

In their study of conformally flat manifolds \cite{SY}, Schoen and Yau showed
that the developing map of any compact, conformally flat manifold 
$M$ admitting a compatible metric of positive scalar curvature, 
is injective as a map from the universal cover of the manifold into 
the $n$-sphere. The corresponding holonomy representation is a monomorphism 
of the fundamental group into the conformal group, $SO(n+1,1)$. 
In particular, the universal 
cover of any such manifold is conformally equivalent to a domain in the 
sphere. They also showed that the complement of the image of this map in
$\SS^n$ has Hausdorff dimension less than or equal to $(n-2)/2$. What they
actually show is that if $g=u^{4/(n-2)}g_{0}$ is a complete metric
of nonnegative scalar curvature on a domain $\SS^n \setminus \Lambda$, 
where $g_0$ is the standard metric on the sphere, and if $g$ has bounded 
Ricci curvature (which is automatically satisfied if $g$ is the lift
of a metric on a compact manifold), then $u\in L^{\frac{n+2}{n-2}}(\SS^{n})$ 
and $u$ extends to a global weak solution of the equation
$$
\Lap_{g_0}u- \frac{n(n-2)}{2}u + \frac{(n-2)}{4(n-1)}
R(g)u^{\frac{n+2}{n-2}} = 0.
$$
Here $R(g) \ge 0$ is the scalar curvature function of $g$. This
leads to the estimate that $\dim \Lambda \le (n-2)/2$. 

Inspired by this result there has been much subsequent work
on case {\bf III}, which is sometimes called the singular Yamabe
problem. Here one is given a compact Riemannian manifold $(M,g_0)$ 
and a closed subset $\Lam \subset M$, and one seeks a metric
$g$ which is conformal to $g_0$, complete on $M \setminus \Lam$
and has constant positive scalar curvature (or CPSC). This
is equivalent to asking for a solution of the equation
$$
\Lap_{g_0}u - \frac{(n-2)}{4(n-1)}R_{0} u + \frac{(n-2)}{4(n-1)}R 
u^{\frac{n+2}{n-2}} = 0 \tag 2.3
$$
on $M \setminus \Lam$, where $u$ is required to blow up sufficiently
strongly on approach to $\Lam$ to ensure completeness of $g$. This last
requirement should be regarded as a sort of singular boundary
condition (albeit at a `boundary' of high codimension). In this
equation $R_0$ is the scalar curvature function of the background metric $g_0$
and $R$ is the prescribed constant scalar curvature of the metric 
$g=u^{4/(n-2)}g_0$. 

It is known that in order for there to exist a solution
to (2.3) with $R \ge 0$, the background metric must satisfy 
an extra condition, namely that the first eigenvalue of
the conformal Laplacian 
$$
 - \Lap + \frac{n-2}{4(n-1)}R(g_0)  \tag 2.4
$$
be nonnegative. This is equivalent to the existence of a metric
(without singularities) of nonnegative scalar curvature
conformal to $g_0$. 

There has also been extensive work concerning this problem when
the the prescribed constant scalar curvature is negative. This
has a quite different character than the positive case, and the
techniques used to produce solutions are also quite different.
For example, one no longer need require any conditions on the
conformal class $[g_0]$; also if a (complete) solution with $R < 0$ exists
then it is necessary that $\dim \Lam > (n-2)/2$. In this case
it is possible to use the maximum principle and barrier arguments, 
and one also finds that solutions are unique. 
Another difference in this negative case is that solutions no longer 
extend to global weak solutions on $M$. A survey of this
case is to be found in \cite{Mc}, and there is some discussion
in \cite{MS}, \cite{M2}. Since then there has been some nice 
recent work by D. Finn \cite{F}. He shows that for any manifold $M$, if
$\Lam$ is a finite, possibly intersecting, union of submanifolds with 
boundary, each with dimension greater than $(n-2)/2$, then 
there exists a complete metric of constant negative scalar
curvature on $M \setminus \Lam$. Recently Finn has also
obtained results concerning existence and nonexistence 
when $\Lam$ is a more general stratified set. We remark
that the issue in this case is always about the 
completeness of the metric $g$, for it is relatively
straightforward to produce solutions to (2.3) with $R < 0$
by barrier arguments. 

When $R = 0$, (2.3) reduces to a linear equation, and the
techniques of \cite{M1} may be used directly, at least when 
$\Lam$ is a (possibly infinite and intersecting) 
union of submanifolds without boundary. It is 
necessary that the dimension of each component of $\Lam$ is
less than or equal to $(n-2)/2$.

Entirely different ideas must be used to produce complete metrics
of constant positive scalar curvature $M \setminus \Lam$. 
In this context the use of gluing techniques to produce new solutions 
began with Schoen's deep construction \cite{S} of 
complete conformally flat metrics on $\SS^{n}\setminus \Lam$,
where $\Lam$ is either an arbitrary finite set, of cardinality
greater than one, or certain nonrectifiable `limit sets'. 
There are many similarities between this and Kapouleas' first
construction of CMC surfaces \cite{K1}, which was developed 
around the same time, although in many technical aspects the two cases
are quite different.
The building blocks here are spheres of radius one, and 
an approximate solution is formed by joining together
infinitely many of these in a certain balanced configuration,
connecting adjacent spheres by very small necks. Just as in the CMC case, 
one must compensate for the fact that these building blocks 
are degenerate in the sense that the linearization of (2.3), $\Cal L$,
has nullspace on each of them. These nullspaces together 
span an infinite dimensional space $X$ such that the restriction
of $\Cal L$ to $X$ has arbitrarily small norm relative to the neck size
if all necks in the approximate solution are sufficiently small. 
When $R$ is normalized to equal $n(n-1)$, this operator takes the form
$$
\Cal L = \Lap + n. \tag 2.5
$$
To perturb these approximate solutions to an exact solution,
Schoen solves the nonlinear problem orthogonally to $X$ and
then shows that if the initial configuration is chosen
correctly this is actually a solution to the full problem. 
More detailed exegeses of this construction are given in
\cite{MPU1, \S 7} and \cite{P1}. A similar approach was used by the 
second author in \cite{P1} to demonstrate nonuniqueness and  
construct high energy solutions for the Yamabe problem on compact 
manifolds.

There are simple examples of CPSC metrics on $\SS^n \setminus \SS^k$:
the standard metric is conformally equivalent, by stereographic projection
from some point on $\SS^k$, to the flat metric on $\RR^n \setminus \RR^k$;
writing this in cylindrical coordinates $(r,\theta,y)$ around $\RR^k$ as 
$dr^2 + r^2d\theta^2 + dy^2$, we conformally change again 
by dividing by $r^2$ to obtain the product metric $r^{-2}(dr^2 + dy^2)
+ d\theta^2$ on $\Bbb H^{k+1} \times \SS^{n-k-1}$. This is a product
of constant curvature spaces, and has constant scalar curvature 
$R = (n-2)(n-2k-2)$ which is positive if and only if $k < (n-2)/2$. 

When $k = 0$ the CPSC metric we obtain this way is the cylinder. Just as 
for CMC surfaces discussed earlier, the cylinder lies in a one-parameter 
family of spherically symmetric metrics, all conformally related
to one another, which are periodic and degenerate to an
infinite set of spheres of radius 1, mutually tangent and
lying along a fixed axis. In analogy to the CMC case we
call these Delaunay metrics $D_\e$, but they correspond
to solutions of an ODE which was first studied by Fowler \cite{Fo1},
\cite{Fo2}. This ODE is rather simple: using the metric
$dt^2 + d\theta^2$ on the cylinder $\RR \times \SS^{n-1}$, then
a conformally related metric $u^{4/(n-2)}(dt^2 + d\theta^2)$ is spherically 
symmetric and has CPSC provided $u$ depends only on $t$ and 
$$
\del_t^2 u - \frac{(n-2)^2}{4} u + \frac{n(n-2)}{4} u^{\frac{n+2}{n-2}} = 0.
\tag 2.6
$$
This equation may be analyzed by writing it as a first order
system for $u$ and $\dot u=\del_t u$. The Hamiltonian energy 
$$
H(u,\dot u) =  \frac{\dot u^2}{2} - \frac{(n-2)^2}{8} u^2  + \frac{(n-2)^2}{8} 
u^{\frac{2n}{n-2}} \tag 2.7
$$
is constant along orbits of this system. The solutions lying within
the set where $H = 0$ are $u \equiv 0$ and $u \equiv \cosh^{(2-n)/2}t$,
the latter of which corresponds to the metric $(\cosh t)^{-2}\, (dt^2 +
d\theta^2)$, which is the incomplete spherical metric on  
$\SS^n \setminus \SS^0$. All solutions
of (2.6) which remain positive for all $t$ remain within the
bounded set $\{ H < 0 \}$, and correspond to periodic orbits,
which are the Delaunay-Fowler solutions. The stationary point in
this set, where $u$ is constant, corresponds to the product metric 
on the cylinder, rescaled to have scalar curvature $n(n-1)$.

One reason for the importance of the Delaunay metrics is the result,
analogous to the result of \cite{KKS} for CMC surfaces, that if
$g$ is a conformally flat metric on a ball with an isolated
singularity at the origin, then $g$ is asymptotic at its
singular point to a unique Delaunay solution $D_{\e}$. This
was proved initially by Caffarelli-Gidas-Spruck \cite{CGS},
and later reproved, with a better estimate of convergence,
by Aviles, Korevaar and Schoen \cite{AKS}. Another proof of 
the Caffarelli-Gidas-Spruck estimate has recently been given
by C.C. Chen and C.S. Lin \cite{CL2}.
The improvement of this estimate to the one obtained by 
Aviles, Korevaar and Schoen is exactly analogous to the 
argument given in \cite{KKS} for case {\bf II}.
It is unknown whether
this result holds for nonconformally flat metrics, although it
is easy to produce, by perturbation methods, nonconformally
flat metrics on a punctured ball with Delaunay asymptotics at the
origin. 

By this asymptotics result we see in particular, that any of the 
solutions produced by Schoen \cite{S1} on $\SS^n \setminus \Lam$, where
$\Lam$ is a finite set, are metrics with asymptotically
Delaunay ends.  The asymptotic Delaunay parameters which arise in
Schoen's construction are all very close to 0. 

The local asymptotics near a higher dimensional singularity
may be much more complicated. It turns out that the product
solutions on $\Bbb H^{k+1} \times \SS^{n-k-1}$ may be extended to
solutions periodic with respect to certain discrete groups
of hyperbolic motions. This is accomplished by regarding this
metric as the universal cover of some compact manifold
$X_\Gamma = (\Bbb H^{k+1}/\Gamma) \times \SS^{n-k-1}$. This manifold
has CPSC, but it is possible to show \cite{M2} that for
certain choices of $\Gamma$ this product hyperbolic-spherical
metric does not minimize the Yamabe energy, hence there is
another CPSC metric which does minimize this energy; its lift
to $\SS^n \setminus \SS^k$ is periodic with respect to the $\Gamma$
action on $\Bbb H^{k+1}$, and it may be shown that it is
constant along the $\SS^{n-k-1}$ factors. 

It is unknown whether these periodic solutions are anomalous,
or whether they exist on more general sets $M \setminus \Lambda$
where $\Lam$ is a $k$-dimensional submanifold. In general
it is possible to use perturbation methods to include a `nondegenerate', 
`admissible' solution metric on $M \setminus \Lambda$ into an
infinite dimensional family (in fact a Banach manifold) of
solutions with the same singular set, where
$\Lambda$ is a nonintersecting union of submanifolds of dimensions
between $0$ and $(n-2)/2$, at least one of which has positive
dimension. Here nondegenerate means that the linearized
scalar curvature operator $\Lap_g + n$ has no $L^2$ nullspace.
Admissible means that the conformal factor $u$ (relative
to the background metric $g_0$ which extends smoothly across $\Lambda$)
has the asymptotic form $u \sim Ar^{(2-n)/2}( 1 + o(1))$, where
$A$ is a purely dimensional constant. These perturbation methods
were first developed for the product metrics on $\SS^n \setminus \SS^k$
in \cite{MS}, but they may also be applied to the solution metrics
constructed in \cite{MPa1}. The $\Gamma$-periodic solutions discussed
above are not admissible, and it would be extremely interesting,
and probably quite difficult, to determine if they are nondegenerate,
hence deformable to a family of asymptotically periodic metrics.

We should mention that if a solution metric $g$ is admissible,
then it is proved in \cite{M2} that $u$ admits a partial
expansion near any component of $\Lam$ with positive dimension $k$ 
of the form
$$ 
u \sim A r^{(2-n)/2}\left(1+r^{k/2 + i\mu} v_+ + r^{k/2 - i\mu}v_-+ 
o(r^{k/2})\right).
$$
The leading asymptotics $v_+$, $v_-$ are functionally related
via a sort of nonlinear scattering matrix; when they are smooth,
there is a full asymptotic expansion, but in general they are only
distributional (of negative order!). Another interesting problem
would be to show that there is a well-posed Dirichlet problem,
i.e. that there is a solution with $v_+$, say, prescribed arbitrarily.
It is possible to show that the set of pairs $(v_+, v_-)$ form
a Lagrangian submanifold of an appropriate infinite dimensional 
function space, so what is required for the Dirichlet problem
is to determine whether this Lagrangian projects bijectively onto the
first `coordinate'. 

When dealing with solutions without such asymptotics one may
seek to determine general upper and lower bounds for solutions, 
in terms of some negative power of the distance to the singular set.
There is a general upper bound for solutions on $\SS^n \setminus \Lam$ of 
the form $u \le C d^{(2-n)/2}$, where $C$ is independent of $\Lam$ and $u$,
and $d$ is the distance function to $\Lam$ \cite{P2}. This does not require 
any regularity of $\Lam$, but it does use the global completeness of the 
metric in a crucial way. Recently C.C. Chen and C.S. Lin \cite{CL1} 
obtained an upper bound of this form for conformally flat solutions 
defined only within a ball. Lower bounds are much more subtle, 
and known in only a few cases. For example, for conformally
flat solutions with isolated singularities, lower
bounds (depending on the solution) are known. Such
solutions are asymptotically Delaunay and there are lower
bounds {\it depending on the necksize} for Delaunay solutions. 
Actually, determining some lower bound for a solution with
isolated singularity is an important step in showing that it is
asymptotically Delaunay. For singular sets $\Lam$ of positive 
dimension, the results of \cite{MPa1} show that there cannot
be lower bounds independent of the solution. 

Recently several existence results have been proved for this problem
using gluing methods. Initially, Pacard \cite{Pa} showed that there are
solutions of this problem on $\SS^n \setminus \Lam$ where $\Lam$
is a finite disjoint union of submanifolds of dimension $(n-2)/2$. 
The paradigm is the usual one: first he constructs a one-parameter
family of approximate solutions which become more and more concentrated
around $\Lam$ and then shows that when these approximate solutions
are sufficiently concentrated, they may be perturbed to exact
solutions. The crucial observation here is that approximate solutions 
may be obtained by pasting radial singular solutions of the equation 
$\Lap u + u^{(n+2)/(n-2)}=0$ on $\RR^N$, $N = n - (n-2)/2 = (n+2)/2$, 
into each fibre of the normal bundle of $\Lam$. The exponent
$(n+2)/(n-2)$ is subcritical in $\RR^N$, and the equation has solutions
 which have damped oscillations as $r \rightarrow \infty$ (as opposed
to the solutions of the `frictionless' critical exponent equation
(2.3) whose oscillations are undamped, i.e. periodic) 
and are concentrated more and more 
at the origin. The error term incurred by pasting these into the 
fibres of the normal bundle of $\Lam$ goes to zero as they become
more concentrated. The importance of the dimension $(n-2)/2$ is that
the equation to be solved is of the form $\Lap_{\SS^n}u + Vu=0$,
where the singularity of $V$ is relatively weak, so that classical
methods may be used to obtain a solution. This method was extended
by Reba\"\i \ \cite{Reb} to consider the case when 
the dimensions of the components of $\Lam$ are allowed to be
slightly less than $(n-2)/2$. These methods are limited by
the fact that only when this dimension is greater than some
$k_0 < (n-2)/2$ do the approximate solutions have finite
instability. More specifically, the linearized scalar curvature 
operator at one of these solutions has only finitely many eigenvalues 
less than zero. By contrast, if the dimension
is less than $k_0$ then the instability becomes infinite;
indeed, in these latter cases $0$ lies in the middle of a band
of continuous spectrum. Nevertheless, replacing the classical
methods in this argument by the pseudodifferential edge theory
developed in \cite{M1}, the first author and Pacard were
able to accomplish the second step for any positive $k$ less than
$(n-2)/2$, perturbing the approximate
solution to an exact solution; in fact, this argument works
for any background manifold $M$ of nonnegative scalar curvature. 
This argument has many delicate
points not encountered in the previous cases, since a subtle
harmonic analysis of the subcritical equation in $\RR^N$ for
radial solutions singular at the origin is now required. 

The result of \cite{MPa1} may be regarded as saying that
the radial singular solutions on $\RR^N$, transplanted
to an $\e$-neighbourhood of the singular set $\Lam$, may
be glued to the zero `solution' (which is of course not even
a metric) on the complement of this tubular neighbourhood. 
The crucial point is that both of these solutions which
are being glued are nondegenerate in the sense that the nullspaces of
their Jacobi operators contain no elements with the appropriate
growth restrictions. This argument was
recently reinterpreted to prove a rather general result \cite{MPU2}
that from any two nondegenerate CPSC metrics, $(M_1,g_1)$ and 
$(M_2,g_2)$, which are either complete or where one or both of the
$M_i$ have boundary, one may construct a (complete) CPSC metric
on the connected sum $M_1 \# M_2$. We wish to discuss this
result a bit further, and in particular to give the correct
general notion of nondegeneracy. 

\head Connected sums and Nondegeneracy \endhead

The Jacobi operator $\Cal L$ 
(with Dirichlet boundary conditions if $\del M \neq 
\emptyset$) of a complete CPSC metric $(M,g)$ is self-adjoint
on $L^2$. Assume that there is some weight function $\a$ on $M$
such that for the corresponding weighted Sobolev spaces 
$$
H^s_{\gamma} = \{ \phi = \a^\gamma \tilde \phi: \tilde \phi \in
H^s \}
$$
one has that 
$$
\Cal L: H^{s+2}_{-\d} \longrightarrow H^s_{-\d}
$$
is injective with closed range for $\d > 0$. Our actual
hypothesis is that for all $s \in \RR$ there is an estimate
$$
\| \phi \|_{s+2,-\d} \le C \| \Cal L \phi \|_{s,-\d} \tag 2.8
$$
for all $\phi \in H^{s+2}_{-\d}(M)$. Note that (2.8) implies
both the closed range and the injectivity.  The dual of $H^{s}_{-\d}$
is $H^{-s}_{\d}$, hence (2.8) also implies that $\Cal L$ is
surjective (and in particular has closed range) on
$H^s_\d$ for all $s \in \RR$. Thus if $f \in H^s_{\d}$ then
there is a solution $v \in H^{s+2}_\d$ to $\Cal L v = f$.
To complete the definition of nondegeneracy 
we also assume that either the nonlinear scalar curvature
operator satisfies
$$
N: H^{s+2}_\d \longrightarrow H^s_\d  \tag 2.9
$$
or else that there exists a `deficiency space' $W$ and a natural 
extension of $N$ so that 
$$
N: H^{s+2}_{-\d} \oplus W \longrightarrow H^s_{-\d} \tag 2.9'
$$
is a well defined, analytic map with surjective linearization
$$
\Cal L: H^{s+2}_{-\d} \oplus W \longrightarrow H^s_{-\d}. \tag 2.10
$$

These nondegeneracy hypotheses (2.8) - (2.10) on the summands $M_1$ and 
$M_2$ are sufficient to ensure that the construction of a metric of 
CPSC on the connected sum $M_1 \# M_2$ may be carried out \cite{MPU2}.
For the specific types of metrics we have discussed above, we remark
that the hypotheses (2.9') and (2.10) are unnecessary for the solution metrics
on $M \setminus \Lam$ from \cite{MPa1} when all components of $\Lam$
have positive dimension.  It is unknown whether the metrics produced
by Schoen \cite{S1} on the complement of a finite set in $\SS^n$
satisfy (2.8) - (2.10). However, the Delaunay solutions do satisfy these
conditions.  To see this, we first define the deficiency space $W$ for
a Delaunay solution $D_\e$.  The elements of $W$ are formed from the
localizations to each end of $D_\e$ of the temperate Jacobi fields.
The Jacobi operator on $D_\e$ takes the form, in cylindrical
coordinates
$$
\Cal L_\e = \Lap_{g_\e} + n = u_\e^{-\frac{2n}{n-2}} \del_t (
u_\e^2 \del t) + u_\e^{-\frac{4}{n-2}} \Lap_\theta + n.
$$
Here $g_\e = u_\e^{4/(n-2)}(dt^2 + d\theta^2)$. Localizing
to each eigenspace of $\Lap_\theta$ leads to an operator with
two solutions $\Phi_{k}^\pm$, indexed by the eigenvalues $\lam_{k}$ 
of $\Lap_\theta$,   one of which (when $k \ne 0$) decays 
exponentially as $t \rightarrow -\infty$ and grows exponentially
as $t \rightarrow \infty$, and the other of which grows and decays,
at the same rate, in the opposite directions. When $k = 0$ the two 
solutions no longer decay or grow exponentially. Indeed,
$$
\Phi_0^+ = \del_t u_\e(t), \qquad \Phi_0^- = \del_\e u_\e(t). \tag 2.11
$$
Thus these Jacobi fields are integrable in the sense that they
correspond to one-parameter families of CPSC metrics on the cylinder:
$\Phi_0^+$ corresponds to the family of `translates' of $g_\e$,
while $\Phi_0^-$ corresponds to the family where the Delaunay
parameter $\e$ varies. The only other Jacobi fields which are
integrable in this sense are the ones corresponding to the first
nonzero eigenvalue on $\SS^{n-1}_\theta$.

For a solution $g$ on $\SS^n \setminus \Lam$, with 
$\Lam = \{p_1, \dots, p_k\}$, a neighbourhood of each singular 
point $p_j$ is asymptotic
to a Delaunay metric $g_{\e_j}$. Thus the truncated Jacobi fields 
$\chi \Phi_0^\pm$, for $\chi$ a cutoff function, on $D_{\e_j}$ 
may be transplanted to a neighbourhood
of $p_j$. The deficiency space $W$ is the $2k$-dimensional linear 
span of these transplanted truncations at all singular points. With
this definition it is proved in \cite{MPU1} that condition (2.10)
holds if the Jacobi operator has a trivial $L^2$ nullspace.
To verify (2.9') we appeal to the fact that the elements
of $W$ are asymptotically integrable, i.e. that up to an
error term in $H^s_{-\d}$ they correspond to one-parameter
families of CPSC metrics on small neighbourhoods of the $p_j$.
Thus, for a small element $w \in W$ we may form a metric
$g_\e(w)$ which has CPSC on the complement of a neighbourhood
of the punctures, and has asymptotically CPSC at each $p_j$,
but the asymptotic Delaunay model at $p_j$ has been translated
slightly and its Delaunay parameter slightly changed according
to the coefficients of $\Phi_0^\pm$ at $p_j$ in $w$. 
With this definition (2.9') now holds.

The general gluing result of \cite{MPU2} states that if 
$(M_1,g_1)$ and $(M_2,g_2)$ satisfy the hypotheses (2.8) - (2.10),
and if $M_\eta \equiv M_1 \#_\eta M_2$ is the connected sum, 
formed by removing
small discs of radius $\sqrt \eta$ around two specified
points $q_j \in M_j$ and identifying the boundaries, then
there is a small perturbation of the resulting approximate
solution metric $g_\eta$ which has CPSC. Furthermore, the
resulting CPSC metric is also nondegenerate.

The main part of the proof proceeds in two steps. Denoting by
$\Cal L_\eta$ the Jacobi operator on the approximate solution 
$g_\eta$, we first show that $\Cal L_\eta$ is
injective on $H^s_{-\d}(M_\eta)$ for $\eta$ sufficiently
small. This means that there is a right inverse (which is
not unique) for $\Cal L_\eta$ as an operator on $H^s_\d$. 
The second step is to show that the norm of an appropriate choice
of right inverse $G_\eta$ does not blow up as $\eta \rightarrow 0$.
The appropriate choice of right inverse is the one with range
lying in the orthogonal complement of the kernel of $\Cal L_\eta$
on $H^s_\d$. This orthogonal complement may be identified
with the range of the adjoint of $\Cal L_\eta$. The way we
prove both of these steps is to assume they are not true
for arbitrarily small $\eta$ and to arrive at a contradiction.
Thus, for step one, we assume that there is a sequence
$\eta_j$ tending to zero and an element $\phi_j \in H^s_{-\d}(M_{\eta_j})$
such that $\Cal L_{\eta_j}\phi_j = 0$. If this were true then we
could obtain, after some normalizations, a nonzero solution
$\phi \in H^s_{-\d}$ on either $M_1$ or $M_2$ satisfying $\Cal L \phi = 0$,
and this contradicts the nondegeneracy hypothesis (2.8) on either
of these summands. For step two, we assume that there are functions
$f_j \in H^s_\d$, $\phi_j \in H^{s+2}_\d$ and $v_j \in H^{s+4}_{\d}$
such that
$$
\Cal L_{\eta_j} \phi_j = f_j, \qquad \Cal L_{\eta_j}^* v_j = \phi_j
$$
where $\Cal L_{\eta_j}^*$ is the adjoint of $\Cal L_{\eta_j}$, 
with $\|f_j\| \rightarrow 0$ and $\| \phi_j \| \equiv 1$. By again
taking limits in a judicious manner, we obtain functions
$v \in H^{s+4}_\d$, $\phi \in H^{s+2}_\d$ satisfying
$\Cal L \phi = 0$, $\Cal L^* v = \phi$ on either $M_1$ or $M_2$, which
is again a contradiction.
A small additional argument is needed to show that the map
(2.10) is surjective on $M_{\eta}$ for $\eta$ sufficiently
small and that it too is `uniformly surjective.'
The remainder of the theorem, the actual perturbation of
the approximate solution to an exact solution, is now straightforward
using (2.10). 

The point we wish to emphasize in this argument is that no explicit
estimates beyond (2.8) are required. Although we do not obtain
explicit bounds for the right inverse $G_\eta$, these are not
really required. This general paradigm, including the indirect
arguments we have sketched, represents a substantial simplification
of the general gluing procedure, at least when the summands
are nondegenerate.

It is not clear whether any of the solutions constructed by
Schoen on $\SS^n \setminus \Lam$ are nondegenerate in this sense. 
More specifically, we have indicated that (2.9') and (2.10) hold for 
these CPSC metrics provided (2.8) is valid. It is shown in \cite{MPU1} 
that the nullspace 
of $\Cal L$ on $H^s_{-\d}(\SS^n \setminus \Lam)$, where $\Lam$ 
is a finite set, is always finite dimensional. Because of 
the parametrices constructed for $\Cal L$ in that paper, to
verify (2.8) it is sufficient to know that this nullspace
is actually trivial.  The only cases where this was known
explicitly were the Delaunay solutions themselves. Thus the
soft gluing theorem of \cite{MPU2} may be applied to show
that $k$ different Delaunay metrics $D_{\e_1}, \dots,
D_{\e_k}$ may be glued together to produce a CPSC metric
on the complement of $2k$ points in $\SS^n$. Moreover, by
being more careful about the definition of the deficiency space 
$W$ in the construction, we show there that it is possible to 
produce  this glued solution in such a way that the Delaunay 
parameters $\e_j$ are unchanged on one
end of each $D_{\e_j}$. This produces the first known
examples of CPSC metrics other than the cylinder with at least some 
cylindrical ends. The analogue of Schoen's problem
mentioned earlier for CMC surfaces, about whether a 
singular Yamabe metric on $\SS^n \setminus \Lam$ with
all ends asymptotically cylindrical must in fact be a cylinder, 
is also open.

This general paradigm for gluing can be applied in some
degenerate problems too. Recently the first author 
and Pacard \cite{MPa2} have been able to reprove Schoen's theorem
concerning solutions with isolated singularities 
and extend it to general manifolds with nonnegative
scalar curvature. Uhlenbeck \cite{U} has announced a proof of the
existence of solutions on $M \setminus \{p\}$ for any
nonnegative scalar curvature manifold $M$ which is not conformally
equivalent to the standard sphere (where no such solution exists). 
The point of both these papers is that while Schoen needs to use
infinitely many balancing conditions in his construction,
essentially one for each small neck connecting the
infinitely many spheres, really only finitely many
balancing conditions are needed. These solutions are
constructed by gluing half-Delaunay metrics onto the
zero solution of $M$. The exact solution will then be
very small, but positive, away from the points $\{p_1, \dots, p_k\}$
and asymptotically Delaunay near each $p_j$. The reason
the degeneracy can be handled is that the elements of the 
relevant deficiency space $W$ (which is slightly larger than
the one we discussed earlier) all correspond to specific
geometric changes of asymptotic geometry. 

\specialhead 3. Moduli spaces \endspecialhead

Once the existence of solutions for any of the problems has been
obtained, it is natural to study the moduli space, or the
set of all solutions. Unlike in the previous section we now 
discuss all three cases together. The moduli space for
solutions of problem {\bf  I} was recently analyzed by 
Ros \cite{Ros} and P\' erez and Ros 
\cite{PR}, while for problems {\bf II} and {\bf III} it was studied 
in \cite{KMP}, \cite{MPU1} and \cite{MPU2}. 
The methods used by P\' erez and Ros were 
slightly different than those used in the other three papers, but it is 
possible to extend these other methods to case {\bf  I} as well.

In what follows we  let $\Cal M_{\Min,g,k}$,
$\Cal M_{\Cmc,k}$ and $\ML$ denote
the moduli spaces of complete, properly immersed minimal surfaces of
finite total curvature with $k$ horizontal ends and genus $g$, 
of complete embedded constant mean 
curvature surfaces with $k$ ends, and of complete metrics of constant positive
scalar curvature on $M \setminus \Lam$ where $\Lam = \{p_1, \dots,
p_k\}$ respectively. In cases {\bf I} and {\bf II} we consider two surfaces
the same if they differ by a rigid motion (thus the dimension of
$\Cal M_{\Min,g,k}$ in (3.1) below is slightly different than in \cite{PR}).
When we discuss results that apply to any of the three
cases, we often write simply $\Cal M_k$ for the relevant
moduli space. (Later we also discuss the `unmarked moduli
space' for case {\bf III}, which will be denoted $\MMk$.) We do not 
discuss the moduli space of solutions for problem {\bf III} when $\Lam$ has
positive dimension. As we have already indicated, in this
case there is an infinite dimensional space of solutions to
the problem. There is hope that this infinite dimensional
moduli space may be parametrized using an asymptotic
Dirichlet problem, but this is as yet unknown. 

\head General features of the moduli spaces \endhead

The basic results proved in these papers are that in each
of these three cases $\Cal M_k$ is locally a real analytic set,
and when $\Sigma \in \Cal M_k$ is a nondegenerate element,
then a neighbourhood $\Cal U$ of $\Sigma$ is a real analytic
manifold. The formal dimensions, which are the actual dimensions
of these neighbourhoods $\Cal U$ around nondegenerate points,
are as follows:
$$
\dim \Cal M_{\Min,g,k} = k, \qquad \dim \Cal M_{\Cmc,k} = 3k - 6, 
\qquad \dim \ML = k. \tag 3.1
$$
To say that $\Cal M_k$ is a locally real analytic set of dimension
$d$ means that near any $\Sigma \in \Cal M_k$ there is a
neighbourhood $\Cal V$ in the space of all nearby surfaces or metrics 
(which are not necessarily solutions) and a real analytic diffeomorphism
$\Psi: \Cal V \longrightarrow \tilde \Cal V$ to an open set
in some finite dimensional space $\RR^m$, such that
$\Psi(\Cal V \cap \Cal M_k)$ is the zero set of a real
analytic function $F: \tilde \Cal V \longrightarrow \RR^\ell$,
where $m - \ell = d$. Thus if the differential of $F$ at
$\Psi(\Sigma)$ is surjective, then $\Psi(\Cal V \cap \Cal M_k)$
is a smooth real analytic manifold of dimension $d$. 

The main point that needs to be understood in each 
case is the possible asymptotic behaviour of solutions.
For case {\bf  I} this is well-known, cf. \cite{HKa}; the ends
of any complete minimal surface in $\RR^3$ are either asymptotically
planar or catenoidal. Moreover, under the requirement that the surface
is embedded, these ends must be parallel.  By only considering
horizontal ends one is simply rotating the surface so that 
the vertical axis is parallel to the axes of these 
asymptotic planes or catenoids.  
The analogous results for complete embedded
CMC surfaces, and for singular Yamabe metrics on $\SS^n \setminus \Lam$,
namely that each end is asymptotically Delaunay, are provided by
\cite{KKS} and \cite{CGS}, \cite{AKS}, \cite{CL2} respectively. 

Because of these results, it is possible to parametrize 
all nearby candidate surfaces or metrics, and to identify
a real analytic mapping of which the solutions, i.e. the
elements of $\Cal M_k$, constitute the zero set. More explicitly,
we use the deficiency space $W$ introduced in the last
subsection for problem {\bf III}, and its analogues for the other cases.
The deficiency space in each case is the set of Jacobi fields
that correspond to one-parameter families of solutions of
the models for each end, suitably truncated and transplanted
to the appropriate end. In each of these problems, these
one-parameter families arise via rigid motions and changes
of the `shape' parameter. 

For problem {\bf  I} the two models are the plane and the catenoid. 
The plane should be thought of as a limiting version of the catenoid,
as the parameter measuring the logarithmic growth of the catenoid
tends to one of its end-values. All catenoids are the same up to
translation, rotation and dilation, so this parameter can be identified
with the diameter of the neck of the catenoid. We shall call it
the logarithmic growth parameter, with the understanding that it
takes the value zero for the plane.  Since we are restricting 
attention to surfaces with horizontal ends, the relevant motions
of each end are translations, either horizontal or vertical, and
change in the logarithmic growth parameter.  By comparing each end to
a fixed horizontal plane in $\RR^3$, we may associate two natural parameters
to each end: the height (e.g. from the central waste of the catenoid
to the fixed plane in the second case) and the logarithmic growth.
We may define the deficiency space $W$ in this case as consisting
of truncations of the asymptotic Jacobi fields corresponding to these changes
of each end; it is $2k$-dimensional. The Jacobi field corresponding
to changing the height is bounded, and in fact tends to the change of 
height, while the Jacobi field corresponding to changing the 
logarithmic growth is linearly growing. By contrast, the asymptotic Jacobi 
fields corresponding to horizontal translations are decaying and are
not included in $W$.  Let $\Cal K_0$ be the space of all Jacobi fields 
which decay along each end. This may be rather large, for example if
there are deformations of the surface where some ends translate horizontally
relative to others. (It is unknown if this phenomenon can occur.)
In any case, $\dim \Cal K_0 \ge 3$, because there are always decaying Jacobi
fields corresponding to the two horizontal translations and to rotation
about the vertical axis. Notice that we do not consider the two other 
asymptotic Jacobi fields on each end, corresponding to rotations not 
preserving the vertical direction, because we are restricting attention
to minimal surfaces with horizontal ends. 

P\' erez and Ros show, using what amounts to
the relative index computation sketched below, that the number $\ell$
of polynomially bounded Jacobi fields on $\Sigma$ is greater than or equal to
$k+3$; in fact, $\ell - (k+3) = \dim \Cal K_0 - 3$. $\Sigma$ is called
nondegenerate when $\ell = k+3$. Hence, when $\Sigma$ is nondegenerate,
all the decaying Jacobi fields correspond to global rigid motions of
$\Sigma$. Even though the Jacobi operator as in (3.2) below is not surjective
when $\Sigma$ is nondegenerate, the obstruction to surjectivity then
consists of the decaying Jacobi fields, which are `geometrically
integrable'. It is not hard to add in these global rigid motions to
make the Jacobi operator surjective, and once this is done,
P\' erez and Ros use the implicit function theorem
to show that $\Cal M_{\Min,g,k}$ is $k+3$-dimensional (or just
$k$-dimensional if surfaces differing by rigid motions are identified)
in a neighbourhood of a nondegenerate surface $\Sigma$.

For problem {\bf II} there are no global decaying Jacobi fields coming
from rigid motions, but there are Jacobi fields, both globally
on $\Sigma$ and for each asymptotic model, arising from the change 
of Delaunay parameter (which is the analogue of scaling the catenoid) 
and from all translations and those rotations orthogonal to the asymptotic 
axes.  Thus the deficiency space $W$ here is a $6k$-dimensional space; for 
each end there is one Delaunay parameter, three translation parameters and 
two rotation parameters.  Now nondegeneracy is the
condition that there are no $L^2$ Jacobi fields at all. Analogous
to the dimension count of \cite{PR}, in this case we prove that
the dimension $\ell$ of polynomially bounded Jacobi fields is
greater than or equal to $3k$. In the nondegenerate case, this
produces a $3k$-dimensional moduli space, although once again
identifying surfaces up to rigid motions we obtain the dimension $3k-6$
as in (3.1). 

Finally, case {\bf III} appears slightly different since it
is intrinsic. The only polynomially bounded Jacobi fields for the asymptotic 
models of each end, i.e. the Delaunay metrics, are the 
Jacobi fields $\Phi_0^\pm$ we have already discussed; these
correspond to changing the Delaunay parameter or translating
along the axis of the Delaunay solution. This translation 
is really with respect to the conformal dilation fixing
the two singular points of the Delaunay metric. Other possible
conformal translations or rotations, because they move the
singular points, correspond to exponentially growing Jacobi fields.
This is in contrast to the other cases above, where these 
geometric Jacobi fields exhibit only linear growth. Thus in this case 
the deficiency space $W$ is a $2k$-dimensional space,
corresponding to the Delaunay and translation parameters.
Below we shall discuss the `unmarked moduli space' for problem
{\bf III}, where the singular points are no longer required
to remain fixed. Then one must consider all Jacobi fields
arising from the conformal action.

To reiterate, the similarities in these three cases is that 
each end has an asymptotic model, and to each model is
associated a set of Jacobi fields which arise geometrically. 
One now builds a finite dimensional space $W$ of asymptotic
Jacobi fields on $\Sigma$, using a partition of unity to 
transplant these model Jacobi fields to each end. $W$ is
called the deficiency subspace. The first main result is that if 
$\Sigma$ is nondegenerate then 
$$
\Cal L: H^{s+2}_{-\d} \oplus  W \longrightarrow H^s_{-\d} \tag 3.2
$$
is surjective (in case {\bf  I} we modify this as discussed above).
The way to prove this is to first observe that
by nondegeneracy and duality, $\Cal L$ is surjective as
a map from $H^{s+2}_\d \rightarrow H^s_\d$. Hence if $f \in
H^s_{-\d}$, then there is an element $u \in H^{s+2}_\d$ such
that $\Cal L u = f$. The second step is to prove a regularity
result for this linear equation, namely that because $f$ decays,
$u$ can be decomposed as a sum of functions $u = v+w$, where
$v \in H^{s+2}_{-\d}$ and $w \in W$. This is proved in \cite{MPU1}
for case {\bf III}, and the same proof works in case {\bf II}. These results
in case {\bf  I} are older, and can be found, for example, in \cite{Mel}.
The Fredholm theory for Laplacians on asymptotically Euclidean or catenoidal
manifolds, or for elliptic operators with the same types of
asymptotic behaviour, dates back to work of Lockhart and McOwen,
and Melrose; the decomposition lemma in this context is found
explicitly in \cite{M1}. 

Now, as described earlier for case {\bf III}, the nonlinear (mean or
scalar curvature) operator $N$ can be defined on $H^{s+2}_{-\d}
\oplus W$ as follows. The elements $w \in W$ of small norm
correspond to nearby surfaces or metrics, $\Sigma(w)$ or $g(w)$,
which are altered on each end according to the components in $w$ of the
truncated Jacobi fields. Thus we obtain these surfaces or metrics
by rotating, translating or changing the relevant parameter (size
of the catenoid or Delaunay parameter) for each end, then
reattaching this slightly altered end to the main body of the
surface in some fixed manner, using cutoff functions. 
The set of all surfaces or metrics which are nearby to $\Sigma$
or $g$ may be written as normal perturbations by some 
decaying function $\phi \in H^{s+2}_{-\d}$ off $\Sigma(w)$
or $g(w)$. Thus we may define 
$$
N(\phi,w) = N_{\Sigma(w)}(\phi), \qquad \text{ or } \qquad 
N(\phi,w) = N_{g(w)}(\phi).  \tag 3.3
$$

The proof that $\Cal M_k$ is a real analytic manifold
of dimension $d$ (one of the three values above) in a 
neighbourhood of a nondegenerate point is now immediate. 
The nonlinear operator in  (3.3) has surjective linearization (3.2), 
and the set of solutions of $N(\phi,w) = 0$ in a small neighbourhood
of zero is a finite dimensional real analytic manifold, by
the implicit function theorem. The dimension of this manifold
is the dimension of the nullspace of the map (3.2). To
determine this we proceed as follows. Because $\Cal L$ is surjective
as a map on $H^{s+2}_\d$ the cokernel is trivial, hence the dimension 
of this nullspace is the same as the index of this map, 
i.e. the dimension of the space of solutions of $\Cal L \phi = 0$ in 
$H^s_\d$ minus the dimension of the cokernel. Even when $\Sigma$ 
or $g$ is not nondegenerate, this index is equal to
the dimension of the `bounded nullspace' $\Cal B$, which we
take as the orthogonal complement of the nullspace of
$\Cal L$ in $H^s_{-\d}$ in the nullspace of 
$\Cal L$ in $H^s_\d$. Now this index could be rather hard to compute, 
but because $\Cal L$ is self-adjoint on $L^2 = H^0_0$, we can identify 
the kernel of $\Cal L$ on $H^s_\d$ with the cokernel of $\Cal L$ on 
$H^s_{-\d}$ and the cokernel on $H^s_\d$ with the kernel on $H^s_{-\d}$.
Hence the index of $\Cal L$ on $H^s_\d$ agrees with one half of the 
`relative index', i.e. the difference between the index on $H^s_\d$ and 
the index on $H^s_{-\d}$. This is fortunate because relative indices are 
much easier to compute 
since they only depend on the asymptotic geometry. In particular,
it is proved in \cite{MPU1} that 
$$
\text{rel-ind\,}(\Cal L)(\d, -\d) = \dim W \qquad \text{ and therefore }
\qquad \dim \Cal B = \frac{1}{2} \dim W.
$$
In case {\bf I}, the bounded nullspace $\Cal B$ is of dimension
$(1/2)(2k) = k$.
In case {\bf II} it gives $(1/2)(6k) = 3k$, and again,
$6$ dimensions are accounted for by rigid motions (since these surfaces
are not required to have horizontal ends as in case {\bf I}), so we
get a dimension count of $3k-6$ for the moduli space.  
Finally, in case {\bf III} it gives
the dimension $(1/2)(2k) = k$.

The proof that in a neighborhood of a degenerate point the moduli
space is still real analytic uses an old method due
to Ljapunov and Schmidt, and Kuranishi. The point is
that (3.1) is not surjective, so we may not apply the implicit
function theorem directly. We may modify (3.1) to get
a surjective linearization as follows. Let $\Cal D$ denote
the finite dimensional `decaying' nullspace of $\Cal L$, i.e. 
the nullspace in $H^s_{-\d}$. Then
$$
\align
&\tilde{\Cal L}: H^{s+2}_{-\d} \oplus W \oplus \Cal D 
\longrightarrow H^s_{-\d}, \\
&\tilde{\Cal L}(v,w,\phi) = \Cal L(v + w) + \phi
\endalign
$$
is surjective. Using this one may produce the real analytic
map $F$ between finite dimensional subspaces. This is discussed
in detail in \cite{KMP}, and the version there is easily adaptable
to cases {\bf  I} and {\bf III}. The earlier proof of this in \cite{MPU1}
for case {\bf III} uses a different sort of argument which doesn't
generalize to the other cases. P\' erez and Ros do not prove
this real analyticity around degenerate points in $\Cal M_{\Min,k}$
explicitly, but note that it follows from the Weierstrass representation.

One crucial point which this general theory does not help address is
whether there are {\it any} nondegenerate points in the moduli space.
Indeed, it is conceivable that $\Cal M_k$ could consist of 
a single degenerate point. It is also not clear whether
there are any degenerate points. It seems very difficult to
give geometric criteria which ensure the nondegeneracy or degeneracy of
a given solution, and it would be very interesting to make
progress on this, or give examples of degenerate or nondegenerate
solutions. In case {\bf I}, P\' erez and Ros note several instances where it
is known that a complete embedded minimal surface is nondegenerate. 
This is because for certain minimal surfaces it is possible
to count the dimension of all Jacobi fields; if this count 
agrees with the formal dimension of the moduli space, then
there can be no decaying Jacobi fields, and hence the surface
is nondegenerate. In the other two cases it is plausible that
any solution that is constructed explicitly, e.g. by gluing,
may be checked explicitly to see if it is nondegenerate.
Unfortunately, it is not at all clear whether the solutions for
case {\bf II} produced by Kapouleas \cite{K1} or for case {\bf III} 
by Schoen \cite{S1} are nondegenerate.  It is hoped that one might
show that solutions for which all necksizes are sufficiently
small (as is true for both Kapouleas' and Schoen's solutions)
are nondegenerate. One scheme for verifying this in case {\bf III} is 
presented at the end of \cite{MPU1}, but whether this scheme works is 
also unclear.  This was one motivation for the construction 
in \cite{MPa2}.  The solutions with isolated singularities constructed 
here are nondegenerate, and this shows that the component of 
$\Cal M_{\Csc, k}$ containing these solutions is of the
predicted dimension. We call these components nondegenerate.
In fact, because this moduli space is
real analytic, it admits a stratification into smooth real
analytic manifolds. The existence of a nondegenerate solution
in a given component means that this solution lies in the
top dimensional `large' stratum of that component. Unfortunately,
it is not necessarily true in the real analytic category
that this component is dense, but at least in the sense of 
measure we can then assert that generic elements in this
component are nondegenerate. 

The solutions in \cite{MPU2} obtained by gluing together
Delaunay metrics are also nondegenerate. Since these are
conformally flat metrics, they may be uniformized and written
as a singular Yamabe metric on some set $\SS^n \setminus \Lam$. 
Not every configuration $\Lam$ can be obtained this
way. Of course, any $\Lam$ here must contain an even number
of elements. More significantly, these configurations depend strongly 
on the sizes of the necks in the connected sum. In fact, as these 
bridging necksizes shrink, these configurations may be written as a union 
of pairs of points. Each pair is widely separated from any other pair,
and the two points within each pair tend toward one another
as the necksizes shrink. For this reason we call these solutions
obtained by gluing together Delaunay metrics `dipole metrics', or
dipole solutions of case {\bf III}. We refer to the singular sets 
$\Lam$ which can be obtained in this way as `dipole configurations'. 
Thus, at least one component of the moduli space $\Cal M_{\Lam}$ 
for any dipole configuration is nondegenerate. 

To use these dipole metrics more effectively, we may 
introduce the unmarked moduli space for case {\bf III}. This space is
actually a closer analogue to the moduli space for case {\bf II} treated 
in \cite{KMP} than the space $\ML$, which would correspond to the submoduli
space of $\Cal M_{\Cmc,k}$ obtained by fixing the asymptotic axes
directions.  The unmarked
moduli space $\MMk$ is defined as the set of all complete conformally flat
CPSC metrics on $\SS^n \setminus \Lam$ where $\Lam$ is
{\it any} configuration of $k$ points on $\SS^n$. Using
modifications of the above arguments, we show in  \cite{MPU2},
that $\MMk$ is a locally real analytic set of dimension $k(n+1)$, 
and there is a natural real
analytic map 
$$
\pi: \MMk \rightarrow \Cal C_k \tag 3.4
$$ 
onto the $kn$-dimensional configuration space of $k$ distinct points on
$\SS^n$; $\Cal C_k$ is naturally identified with an open subset
of $(\SS^n)^k$, and hence is a real analytic manifold. 
We may once again talk about nondegenerate components of 
$\MMk$, and we see from the discussion above that
any component of this unmarked moduli space which contains
a dipole metric is nondegenerate. In particular, we obtain
nondegeneracy of metrics for configurations $\Lam$ which 
are not necessarily dipole configurations. 

The space $\MMk$ is a natural setting for a more global 
investigation of the geometry of these moduli spaces.  This is the subject
of recent work of the second author \cite{P3} which we 
describe below.  We first recall the natural invariants associated with 
the elements of  $\ML$  and the infinitesimal Lagrangian 
structure which these moduli spaces carry as a direct consequence 
of the linear analysis.

\head Poho\v zaev invariants \endhead 

There are some interesting invariants associated to a metric $g \in \ML$
at each of the singular points $p_i\in\Lam$.  These are derived from 
the generalized Poho\v zaev identity derived by Schoen \cite{S1}
which takes the form
$$
\int_{\Omega}L_{X}R(g)\,dv_{g}=
\frac{2n}{n-2}\int_{\partial \Omega}(Ric_{g}-n^{-1}R(g)g)(X,\nu)\,d\s_{g},
$$
where $\Omega$ is a domain with smooth boundary $\partial\Omega$, $X$ is a 
conformal Killing vector field on $\Omega$, $L_{X}$ denotes the Lie 
derivative, and $\nu$ is the outward unit normal vector to 
$\partial \Omega$.  When $X$ is not conformal Killing one has an 
additional interior term.  This identity is a Riemannian 
version of Poho\v zaev's original identity \cite{Po}, which in turn is 
a nonlinear version of the classic Rellich identity \cite{Rel}.  
These are all derived from the divergence theorem. 

Since any $g\in \ML$ has constant scalar curvature $R(g)=n(n-1)$, 
this identity implies that the expression
$$
\frac{2n}{n-2} \int_{\Sigma} ( Ric_g - n^{-1} R(g)g)(X,\nu)\,d\sigma_g
$$
only depends on the homology class of $\Sigma$ in 
$H_{n-1}(\SS^{n}\setminus\Lam;\RR)$.  Fix submanifolds $\Sig_i$ which
are the boundaries of small balls around each singular point
$p_i \in \Lam$. Their homology classes $[\Sig_i]$ constitute a
basis for this space. Representing the space of conformal Killing fields as
the Lie algebra of the conformal group $SO(n+1,1)$, we obtain elements
of the dual Lie algebra defined by 
$$
{\Cal P}_{i}(g)=\int_{\Sig_i}(Ric_{g}-(n-1)g)(\, \cdot\, ,\nu)d\s_{g}
\in\frak{so}^{*}(n+1,1).
$$
These invariants play a central role in the natural compactification
of the moduli space described below.  It is sometimes sufficient to
work with the simpler {\it dilational Poho\v zaev invariants} 
defined by
$$
{\Cal P}_{i}^{0}(g)={\Cal P}_{i}(g)(X_{p_{i}})\in \RR,
$$
where $X_{p_{i}}\in \frak{so}(n+1,1)$ is the generator for the one parameter
family of centered dilations fixing $p_i$ and its antipodal point.  This 
component of the Poho\v zaev invariant may be explicitly computed 
from the asymptotics.  This is done in \cite{P2} and one finds that 
$$
{\Cal P}_{i}^{0}(g)=c(n)H(g_{\e_i})
$$
where $c(n)$ is an explicit positive constant and $H(g_{\e_i})<0$ is the 
Hamiltonian energy of the asymptotic Delaunay metric, as in (2.7).

The advantage of working in $\MMk$ is that there is a natural action of the
conformal group on this space, $SO(n+1,1):\MMk\rightarrow\MMk$, via
$$
(g,\Lam)\mapsto (G^{*}(g),G^{-1}(\Lam)),
$$
for $G\in SO(n+1,1)$.  Under this action the Poho\v zaev invariants
transform under the co-adjoint representation of $SO(n+1,1)$.

There are analogous invariants in case {\bf II}. 
Suppose $\Sigma\in \Cal M_{\Cmc,k}$ and let $\Gamma\subset\Sigma$ be a 
smooth $1$-cycle, and $K\subset\RR^2$ a smooth 
$2$-chain such that $\partial K =
\Gamma$. If $\nu$ is a unit normal to $K$ and $\eta$ is a unit conormal to
$\Gamma$ on $\Sigma$ so that all corresponding orientations are coherent 
with that of $\RR^3$, then for any Killing vector field $X$ in $\RR^3$, 
the expression
$$
\mu([\Gamma])(X)=\int_{\Gamma} \eta\cdot X \,ds + \int_{K} \nu\cdot X\,d\sigma
$$
depends only on the the homology class of $\Gamma$ in $\Sigma$.
This is again a clever application of the divergence theorem coupled
with the first variation formula for $\Sigma$, c.f. \cite{KKS}.  Again, by 
choosing a basis for the first homology of $\Sigma$ one has an invariant 
$\mu_{0}\in \frak{so}^{*}(3)$.  It is easy to see that this invariant 
also transforms under the action of $SO(3)$ by the co-adjoint representation.
This invariant is exploited in \cite{KK} to obtain a priori area and curvature
bounds for $\Sigma\in \Cal M_{\Cmc,k}$.

\head Coordinates on the moduli space \endhead

Here we describe one way of obtaining infinitesimal coordinates on the
moduli space and recall from  \cite{MPU1} and \cite{KMP} how these give rise
to an infinitesimal Lagrangian structure. For simplicity we describe the 
results here for $\ML$, though there are exact analogues for $\Cal M_{\Cmc, k}$
and $\MMk$.  We also describe the Lagrangian structure of P\' erez and Ros 
for $\Cal M_{\Min, k}$.

Suppose $g\in\ML$ is nondegenerate, so that in particular (2.10) holds.
The tangent space $T_{g}\ML$ is identified with the kernel of this
map which we have called the `bounded nullspace' 
$\Cal B \subset H^{s+2}_{-\d} \oplus W$.
Thus any $\phi\in T_{g}\ML$ can be identified (up to an exponentially 
decaying error) with a linear combination of elements of $W$;
$\phi \sim \sum_{j=1}^{k}\left(a_{j}(\Phi_0^+)^{(j)} +  
b_{j}(\Phi_0^-)^{(j)}\right)$.  This determines a map 
$$
\gathered S:  T_{g}\ML\longrightarrow \RR^{2k} \\
 \phi \longmapsto (a_1, b_1, \dots , a_k, b_k). \endgathered 
$$

The assumption that $g$ is nondegenerate implies that $S$ is injective:
for, if $S(\phi)=0$ then $\phi\in H^{s+2}_{-\d}$, 
which would imply that $\phi=0$, since ${\Cal L}\phi=0$ and $g$ 
satisfies (2.8). Hence in this case the image of $S$ is isomorphic 
to the nullspace of $\Cal L$, which is identified with 
$T_g \ML$. Using Green's theorem \cite{MPU1} one shows that this image 
is a Lagrangian subspace of $\RR^{2k}$ endowed 
with the standard symplectic form $\sum_{j=1}^k da_j \wedge db_j$.
Hence a neighbourhood $\Cal U$ of $g$ in $\ML$ inherits these
`Lagrangian coordinates'. This formalism works nearly identically for
$\Cal M_{\Cmc, k}$ and $\MMk$. However, it is not obvious how to
extend this infinitesimal Lagrangian structure to a global one.  This is
discussed in the next subsection.

As we described at the beginning of this section the geometric
parameters for elements of $\Cal M_{\Min, k}$ are the logarithmic 
growth and height parameters.  For $\Sigma\in\Cal M_{\Min, k}$ 
we let $(\text{\bf L}(\Sigma), \text{\bf H}(\Sigma))$ denote this $2k$-tuple.
There is a natural map $f:\Cal M_{\Min, k}\rightarrow\RR^{2k}$ given 
by $\Sigma\mapsto (\text{\bf L}(\Sigma), \text{\bf H}(\Sigma))$.
One may again endow $\RR^{2k}$ with the standard symplectic structure 
as above, so that the logarithmic growth and height parameters correspond
to the numbers $a$ and $b$, respectively.
P\' erez and Ros show that the restriction of this map $f$
to the nondegenerate surfaces is an analytic Lagrangian immersion.
The fact that the map is isotropic follows from Green's theorem. 

\head Global features of the moduli space \endhead 

As a means towards understanding the set of all solutions to these
problems more thoroughly, we would like to explore which natural geometric 
structures, if any, these moduli spaces carry.  One immediate remark is
that these spaces do not support a natural $L^{2}$-metric because the 
Jacobi fields are not square integrable.  Thus these spaces do not have
the natural Riemannnian structure which is common to many other geometric 
moduli spaces.  Another geometric 
structure often found on moduli spaces is a realization of the space as a 
Lagrangian submanifold in a natural symplectic configuration space.
The existence of such a structure for cases {\bf II} and {\bf III} 
is hinted at in the infinitesimal Lagrangian structure described above.
In this subsection we describe some recent work 
\cite{P3} which attempts to extend this to a global picture.  
P\' erez and Ros also exhibit such a structure on $\Cal M_{\Min, k}$.

The relevant configuration spaces for cases {\bf II} and {\bf III} 
are simply the spaces of Delaunay surfaces or metrics.  
For simplicity we discuss only case {\bf III} here.  
In this context, we regard a Delaunay metric as a complete, conformally 
flat CPSC metric on the complement of two points in $\SS^{n}$.  
Any such metric is conformally  equivalent to one of the ODE solutions 
discussed above.  These ODE solutions may be parametrized by a 
noncompact surface $\Omega$ with coordinates corresponding to the 
Delaunay and translation parameters, and corresponding area form 
$\omega_{0} =  da \wedge db$ as above.  To specify a Delaunay metric 
we need to choose two distinct points in $\SS^n$ and an ODE solution.  
Thus the space of Delaunay metrics is simply
$$
{\Cal D} = \left[(\SS^{n}\times \SS^{n})\setminus{\Lap}\right]\times 
\Omega,
$$
where $\Lap$ denotes the diagonal.  The space of pairs of distinct
points in the $n$-sphere, $(\SS^{n}\times \SS^{n})\setminus{\Lap}$, has a 
natural symplectic structure when regarded as the space of geodesics in 
hyperbolic $(n+1)$-space, $\HH^{n+1}$.  The space of $k$ Delaunay metrics 
is thus the $k$-fold Cartesian power
$$
{\Cal D}^{k} = \prod_{i=1}^{k}
\left[(\SS^{n}\times \SS^{n})\setminus{\Lap}\right]_{i}\times 
\Omega_{i}.
$$
${\Cal D}^{k}$ inherits a natural symplectic structure from the product
structure, $\omega$ on ${\Cal D}$, coming from the symplectic structure on 
$(\SS^{n}\times \SS^{n})\setminus{\Lap}$ and the area form 
$\omega_0$ on $\Omega$.
We let $\omega^k$ denote this symplectic form on ${\Cal D}^{k}$. 

There are two natural group actions on ${\Cal D}^{k}$.  The first is an
action of the conformal group $S0(n+1,1)$ acting as hyperbolic isometries 
on the space of geodesics and as, suitably interpreted, conformal
transformations on $\Omega$.  The second is the obvious action of the 
symmetric group on $k$ letters,  {\bf  S}$_k$, permuting the factors
of ${\Cal D}^{k}$.  Each of these groups act as symplectomorphisms on
${\Cal D}^{k}$.

The space ${\Cal D}^{k}$ is the natural configuration space for the 
unmarked moduli space $\MMk$.  In particular
\proclaim{Theorem \cite{P3}} There is a natural $S0(n+1,1)\times 
\text {\bf  S}_k$ equivariant (possibly singular) immersion
$$
\Psi:\MMk\longrightarrow {\Cal D}^{k}.
$$
Moreover, the map $\Psi$ is well defined on all of $\MMk$ and nonsingular 
off of the singular locus.
\endproclaim

This map is constructed using the Poho\v zaev invariants discussed above.
We expect that an equivalent construction can be given analytically by 
using a more refined asymptotic expansion for $(g,\Lam)\in\MMk$.

One should note that the dimension of ${\Cal D}^{k}$ is $2k(n+1)$,
which is twice the dimension $\MMk$.  The claim is that 
$\Psi$ is an equivariant Lagrangian immersion.    
Recall that for $\Lam\in{\Cal C}_k$,
$\ML = \pi^{-1}(\Lam)$, where $\pi :\MMk\rightarrow {\Cal C}_k$ is the 
projection.  Let $i:\ML\hookrightarrow\MMk$ be the inclusion.  The map 
$\Psi$ immediately provides a generalization of the infinitesimal 
Lagrangian structure on $T_{g}\ML$ described above, namely
\proclaim{Corollary} The map 
$\psi=(\Psi\circ i):\ML\rightarrow {\Cal D}^{k}$
is isotropic, i.e. $(\psi)^{*}(\omega^{k})=0$.
\endproclaim

The expectation is that the equivariant Lagrangian immersion of $\MMk$
will provide additional tools to address global questions about the 
moduli space, in particular the structure of its singular locus.
Whether this map is actually an embedding is a difficult open question.
This is related to the question of the uniqueness of elements of $\MMk$ 
with specified asymptotics.  There is an analogous picture 
for case {\bf II} which will be explored elsewhere. As pointed out earlier, 
P\'erez and Ros define a Lagrangian immersion of $\Cal M_{\Min, k}$ into 
$\RR^{2k}$, the space of Log and Height parameters with standard 
symplectic structure, in \cite{PR}. They also compute the second
fundamental form of the image of this immersion. 

\head Compactifications of the moduli spaces \endhead

We have already alluded on several occasions to ways in
which sequences of elements in either $\Cal M_{\Cmc, k}$ or 
$\ML$ can degenerate. These examples indicate that
these moduli spaces are not compact. There are also good
models for degenerations of sequences of elements in 
$\Cal M_{\Min, k}$. We are thus led to the issue of whether any of
these spaces admit natural compactifications. We have studied
this question in some detail for the spaces $\ML$
and we describe these results below. Similar results should 
hold for the moduli spaces of CMC surfaces; this has been partially
explored in \cite{KK}, but we leave the development and application of 
these for elsewhere.  We first briefly describe the results 
of Ros \cite{Ros} in this direction for $\Cal M_{\Min, k}$.

Ros examines two types of compactness. The first he calls weak compactness,
and it really describes a compactification of the moduli spaces
$\Cal M_{\Min,g,k}$. He shows that if $\Sigma_j$ is a sequence
of elements in one of these moduli spaces, then a subsequence
converges to a finite union $\Sigma^{(1)}, \dots, \Sigma^{(\ell)}$
of minimal surfaces with lower genus or fewer number of ends. 
He also shows that some of his moduli spaces are `strongly compact',
i.e. that any sequence $\Sigma_j$ as above has a subsequence converging
to a minimal surface $\Sigma^{(0)}$ in the same moduli space. The criterion to
show that some $\Cal M_{\Min,g,k}$ is strongly compact is simply
that all moduli spaces of surfaces with lower genus or fewer number
of ends are empty. There are several known `nonexistence' results
of this type known. In any case, this strong compactness is a simple
corollary of the more general weak compactness result.

The results for cases {\bf II} and {\bf III} are quite similar
in form, although we never obtain actual compactness of any of
our moduli spaces, simply because it is known that all the
moduli spaces are nonempty. Anyway, we proceed to describe
our results, concentrating on case {\bf III}.

The first basic result gives a geometric
criterion for the convergence of a sequence of elements
$\{g_j\} \subset \ML$.
\proclaim{Theorem \cite{P2}} Suppose that $\{g_j\} \subset \ML$ is
a sequence of complete CPSC metrics on $\SS^n \setminus \Lam$.
Suppose also that the dilational Poho\v zaev invariants at
each singular point $p_i$, ${\Cal P}_{i}^{0}(g_j)$ are bounded
away from zero. Then  a subsequence of the $g_j$'s converge uniformly 
in $C^\infty$ on compact subsets of $\SS^n \setminus \Lam$ to an 
element $g_\infty \in\ML$. 
\endproclaim

A corollary of this result is that the only way in which the
sequence $\{g_j\}$ can degenerate to the `boundary' of the moduli
space is if at some nontrivial subset $\Lam' \subset \Lam$ the necksizes
of these metrics tend to zero. Suppose then that this occurs.
The results of \cite{P2} still imply convergence of the metric
as a symmetric $2$-tensor away from $\Lam$. If $\Lam' \neq \Lam$, 
then the limiting tensor $g_\infty$ is a nondegenerate metric away 
from $\Lam$. However, one can show that $g_\infty$ is smooth
and nondegenerate across the points of $\Lam '$ ; hence it is a
singular Yamabe metric on $\SS^n \setminus (\Lam'')$, where 
$\Lam''=\Lam \setminus \Lam'$.  Thus $g_\infty$ is  an element of 
$\Cal M_{\Lam''}$. When $\Lam' = \Lam$ then one of two
possibilities can occur: either $g_\infty$ is isometric to the standard
metric on $\SS^n$, or else it vanishes identically. An illustration
of this last case is provided by the explicit picture available
for $\ML$ when $\Lam=\{p,q\}$ consists of two points. 
This space can be identified with the interior of the
set $\{H < 0\}$, where the Hamiltonian energy is negative.
The natural compactification of this set is obtained by adding
the boundary $\{H = 0\}$, which is the union of two orbits:
the zero orbit and the one corresponding to the standard metric
on the sphere. So we have actually added on infinitely many
copies of the spherical metric, one for each conformal dilation
fixing the two singular points, but only one copy of the zero solution.
The compactified moduli space $\{H \le 0\}$
is a real subanalytic set. 

Since we have given a description of all possible degenerations of
sequences in $\ML$, we have in fact shown that
there is a compactification $\CML$ obtained
by adding on certain components of the moduli spaces for fewer
numbers of singular points. Exactly which of these `submoduli
spaces' occur in the compactification seems to be  a very 
subtle question. In particular, it is not clear what the Hausdorff
codimension of the full boundary $\CML \setminus
\ML$ is; when $k = 2$ this codimension is only one.
It is also an interesting question whether these compactified
moduli spaces are real subanalytic sets in general. It seems
likely that this might be resolved via the Lagrangian embedding
picture described above.

We may also consider the compactifications of the unmarked
moduli spaces $\MMk$.  As before, these compactifications
may be analyzed by studying sequences of metrics in $\MMk$.
The first result, also proved in \cite{P2}, is
\proclaim{Proposition 3.5} Let $\{g_j\}$ be a sequence of elements in $\MMk$
such that all dilational Poho\v zaev invariants at the $k$ varying
singular points $p_i^j = p_i(g_j)$ are bounded away from zero, and the 
distances, $d_j=\min_{i\neq k}(\text{dist}(p_i^j,p_k^j))$, 
between these singular points are uniformly bounded
away from zero. Then a subsequences of the $g_j$'s converges to 
an element in $\MMk$.
\endproclaim
There are now more possibilities to consider when 
analyzing the possible modes of degeneration for sequences $\{g_j\} 
\subset \MMk$. The first is that some of the dilational Poho\v zaev invariants
tend to zero, but the singular points stay uniformly
bounded away from each other. The analysis of this case is
exactly the same as when the singular points are fixed.
The new possibility is when some of the singular points
tend to one another, and this divides into subcases, 
according to whether some of the dilational Poho\v zaev 
invariants tend to zero or not, and whether the singular 
points where these necksizes are shrinking to zero are
amongst the ones colliding or not. 

Suppose first that some collection of the points are coming together, but
all dilational Poho\v zaev invariants are bounded away from
zero. To be specific, suppose that the points $p_1, \dots, p_\ell$
stay uniformly bounded away from each other, as well
as from the remaining points $p_{\ell+1}, \dots, p_k$, and suppose
that these remaining points converge to some point $q$. 

We first note that we may always assume that $\ell > 1$. To arrange
this, pull back the configuration by a conformal transformation so that 
$p_1$ and $p_2$ are antipodal. If both of these points remain 
a certain distance from the others, then $\ell \ge 2$. If not, then 
there is a cluster converging to $p_1$ say. We can pull this cluster
apart by a further conformal transformation, positioning $p_1$
and some other $p_i$ antipodally. This process must terminate 
after finitely many steps, and by a relabeling of the points
we are in the situation above. 

Since the dilational Poho\v zaev invariants at all the points remain
bounded away from zero, we may extract a convergent subsequence of the
metrics which tend to a limiting smooth metric $g$ away from $\{p_1, \dots, 
p_\ell, q\}$. Because $g$ has isolated singular points, at each point
it is either asymptotically Delaunay or the singularity is removable, 
i.e. the metric extends smoothly across the point.
The former occurs at the $p_i$ because the dilational Poho\v zaev
invariants did not tend to zero, while at $q$ either case may occur.
If $g$ extends smoothly across $q$, we may picture the degeneration 
of the whole sequence as a `bubbling off' of a singular Yamabe metric 
with $\ell$ singularities from a possibly degenerating sequence of singular 
Yamabe metrics with $k-\ell$ singularities. If this sequence with $k-\ell$
singularities does not itself degenerate, then the whole degeneration
corresponds to the shrinking of a neck in a connected sum between two
singular Yamabe metrics, one with $\ell$ singularities and one with
$k - \ell$ singularities. However, even if the solution with $k-\ell$
singularities degenerates, because all dilational Poho\v zaev 
invariants are bounded away from zero, we can repeat this analysis.
The other case, when $g$ is singular at $q$, may be pictured as 
the degeneration of a sequence of connected sums between two metrics
$g_\alpha$ and $g_\beta$, where $g_\alpha$ converges to an element
of $\Bbb M_{\ell+1}$, but $g_\beta$ may possibly degenerate further,
and where these metrics are connected at points $x_\alpha$, $x_\beta$,
where $x_\alpha$ converges to one of the singularities of the
limit of $g_\alpha$. In any event, we may still isolate the 
$k-\ell$ singular points by pulling back with a conformal transformation
and continue the analysis. In the end, we obtain the

\proclaim{Proposition 3.6} Suppose that $\{g_j\}$ is any sequence of elements
in $\MMk$, such that the dilational Poho\v zaev invariants at the singular 
points $p_i^j$ are all bounded away from zero. Then there is a partition
of $\{1, \dots, k\}$ into subsets $I_1, \dots, I_s$, with each $|I_\ell| 
\ge 2$, such that for each $I_\ell$ there is a a sequence of conformal
transformations $f_j$ such that a subsequence of 
$f_j^*(g_j)$ converges to an element of $\Bbb M_{|I_\ell|}$.

\endproclaim
Thus the general picture in this case is of $s$ different singular 
Yamabe metrics bridged together by connected sums, with all connecting
necks shrinking to zero in the limit. The construction in \cite{MPU2}
can be used to give examples of this multi-connected sum degeneration.

When some of the dilational Poho\v zaev invariants tends to zero, it may 
seem that there are further cases to consider. However, a moment's thought
shows that, using the results of \cite{P2}, the only difference
in this case is that in the preceding proposition, having
chosen $I_\ell$ and the conformal transformations $f_j$, the limits of
the metrics $f_j^*g_j$ converge to elements in some unmarked moduli
space $\Bbb M_{m}$ where $m \le |I_\ell|$. Thus the general
statement of this proposition may be altered to cover all cases
by simply requiring $\{I_1, \dots, I_s\}$ to be a partition of
the subset $\{p_{i_1,}, \dots, p_{i_r}\}$ of points where the
dilational Poho\v zaev invariants do not tend to zero.

This completes the analysis of the possible degenerations of sequences
in $\MMk$. The conclusion is that a compactification $\bar{\Bbb M}_k$ may
be obtained by adding on a union of certain subsets of unmarked moduli
spaces of metrics singular at fewer than $k$ points. 

\vfill\eject

\enddocument

\bye